\documentclass[english,9pt]{article}
\usepackage{amsmath,amssymb,amsthm,babel}
\newtheorem{prop}{Proposition}
\begin{document}

\def\proof{{\it Proof.}\ }
\def\pt#1{{\it Proof of Theorem \ref{#1}.}}
\def\eq#1{(\ref{#1})}
\def\th#1{Theorem \ref{#1}}
\def\neweq#1{\begin{equation}\label{#1}}
\def\endeq{\end{equation}}
\def\endproof{\nolinebreak\hfill\rule{2mm}{2mm}}
\def\ri{\rightarrow}
\def\di{\displaystyle}
\def\RR{{\mathbb R} }
\def\NN{{\mathbb N} }
\def\ZZ{{\mathbb Z}}
\def\QQ{{\mathbb Q}}
\def\ep{\varepsilon}
\def\dz{\sharp}
\def\we{\rightharpoonup}
\def\ii{{\^{\i}}}
\def\s{\c s}
\def\ct{\c t}
\def\a{\u a}
\def\aa{\u a\;}
\newtheorem{teo}{Theorem}[section]
\newtheorem{lem}{Lemma}[section]
\newtheorem{corollary}{Corollary}[section]
\newtheorem{rem}{Remark}[section]
\renewcommand{\theequation}{\thesection.\arabic{equation}}
\title{Preliminary results on the homogenization of thin piezoelectric perforated  shells\footnote{This preprint is translation of Chapter 3 of PhD thesis \cite{Mech2} (Full text: http://tel.archives-ouvertes.fr/tel-00008496/fr/)}}

\date{}
\author{Houari Mechkour\thanks{%
Centre de Math\'ematiques Appliqu\'ees (UMR 7641)
\'Ecole Polytechnique, 91128 Palaiseau, France. (mechkour@cmap.polytechnique.fr).}}
\maketitle

\begin{abstract}
We consider a composite piezoelectric material whose reference configuration is a thin shell with fixed thickness. In this work, we give a new approach based on the periodic unfolding method to justify the modelling of a thin piezoelectric perforated  shells and we establish the limit constitutive law by letting the size of holes is supposed to go to zero. This allows to use the homogenization technique to derive the limitting equations and the homogenizaed coefficients are explicity described.\\

\end{abstract}

\section{Introduction}

The shells is a three dimensional continous medium, where its thickness is small compared to other dimensions. Its geometry is characterized by two small parameters : the tickness of the shells and the size of perforations. The behavior of piezoelectric shells when the tickness goes to zero has been studied by Haenel \cite {Haen}.

In this work, we consider a periodically perforated piezoelectric shells. For two dimensional limit equations obtained by Haenel \cite{Haen}, we study the behavior for the elastic displacement and electric potential 
as the size of perforations becomes smaller and smaller.We are concerned with two independent problems : the membrane and bending problems and we give the convergence results based on new {\it periodic unfolding method}, recently introduced by Cioranescu, Damlamian and Griso \cite{Cior}. 

We find explicitly the overall homogenized tensors and we study their properties in order to give
a corrector results associated for membrane and bending problems.
\section{The periodic unfolding method}

The {\it Periodic Unfolding Method} is a novel approach to periodic homogenization 
problems that applies as well to problems with holes and truss-like structures or 
in linearized elasticity. The periodic unfoding method is equivalent to two-scale
convergence method but is both simpler and more efficient.

First we briefly introduce this method which was developped in 
\cite{Cior}.
\medskip

Let $\Omega\subset\RR^3$ be a domain with a smooth boundary
$\partial\Omega$ and let $Y=[0,1]^3$ be
the reference cell. We consider that $S$
is an open smooth boundary subset of  $Y$ such that 
$\overline S\subset \overline Y$ and set
$Y^*=Y\setminus S.$ 
$S$ plays of the reference hole while $Y^*$ is the part of $Y$ 
occupied by the material. We also set 
$S^\ep=\ep(\overline S+k)\cap \Omega,$ $k\in \ZZ^3.$

The perforated domain $\Omega_\ep$ is defined as the set
$\Omega^\ep=\Omega\setminus \overline{S^\ep}.$
We assume that $\Omega^\ep$ is connected and that the holes do not 
intersect the boundary $\partial\Omega.$
\smallskip

 For $z\in \RR^3,$ we denote by $[z]_{Y^*}$ the unique 
integer combination  such that $z-[z]_{Y^*}\in Y^*$ and set
$\{z\}_{Y^*}=z-[z]_{Y^*}\in Y^*.$\\
For any $x\in\RR^3$ and $\ep>0$ we have
$$\di x=\ep\left(\left[\frac{x}{\ep}\right]_{Y^*}+
\left\{\frac{x}{\ep}\right\}_{Y^*}\right)$$
Define ${\cal T}^{\ep}:L^2(\Omega^\ep)\ri L^2(\Omega^\ep\times Y^*)$ 
with
$$\di {\cal T}^{\ep}(w)(x,y)=w\left(\ep\left[\frac{x}{\ep}\right]_{Y^*}+\ep y\right),
\quad\mbox{ for all }(x,y)\in\Omega^\ep\times Y^*.$$

Obviously, for any $v,w\in L^2(\Omega^\ep)$ we have
\neweq{aunu}
{\cal T}^{\ep}(vw)={\cal T}^{\ep}(v){\cal T}^{\ep}(w)
\endeq
\neweq{adoi}
{\cal T}^{\ep}(v+w)={\cal T}^{\ep}(v)+{\cal T}^{\ep}(w)
\endeq
For our purpose, all functions defined in $L^2(\Omega^\ep)$ are extended by zero outside $\Omega^\ep.$
\begin{prop}\label{p1}{\rm (Properties of ${\cal T}^{\ep}$)}{\rm (see \cite{Cior})}
\begin{description}
\item[(a)] For all $w\in L^1(\Omega^\ep)$ we have
\neweq{atrei}
\di \int\limits_{\Omega^\ep}w dx=\frac{1}{|Y^*|}\int
\limits_{\Omega^\ep\times Y^*}{\cal T}^{\ep}(w)dxdy;
\endeq
\item[(b)] For any $w\in L^2(\Omega^\ep)$ we have
\neweq{apatru}
{\cal T}^{\ep}(w)\ri w\qquad\mbox{ strongly in }\quad
L^2(\Omega\times Y^*);
\endeq
\item[(c)] If $(w^{\ep})\subset L^2(\Omega^\ep),$ then

$\di w^{\ep}\we w\;\;\mbox{weakly in }L^2(\Omega)\Longrightarrow
{\cal T}^{\ep}(w^{\ep})\we w\;\;\mbox{weakly in }
L^2(\Omega\times Y^*);$

$\di{\cal T}^{\ep}(w^{\ep})\we \widehat w\;\;\mbox{weakly in }L^2(\Omega\times Y^*)
\Longrightarrow
w^{\ep}\we \frac{1}{|Y^*|}\int\limits_{Y^*}\widehat w dy\;\;\mbox{weakly in } L^2(\Omega).$
\end{description}
\end{prop}
\begin{prop}\label{p2}{\rm (see \cite{Cior})}
Let $(w^{\ep})\subset L^2(\Omega^\ep)$ be a bounded sequence. Then
$$\di {\cal T}^{\ep}(w^{\ep})\we w\;\;\mbox{weakly in}\,L^2(\Omega\times Y^*)
\Longleftrightarrow\;w^{\ep}\;\mbox{ two-scales converges to } w.$$
\end{prop}

\begin{teo}\label{t1}{\rm (see \cite{Cior})}
Let $(w^{\ep})\subset H^1(\Omega^\ep)$ be a bounded sequence that weakly converges to $w$
in $H^1(\Omega).$ Then there exists $ w^1\in L^2(\Omega; H^1_{\rm per}(Y))$ such that
\neweq{acinci}
{\cal T}^{\ep}(w^{\ep})\we w\quad\mbox{ weakly in }\,L^2(\Omega\times Y^*),
\endeq
\neweq{asase}
{\cal T}^{\ep}(\nabla_x w^{\ep})\we \nabla_x w+
\nabla_y w^1
\quad\mbox{ weakly in }\,L^2(\Omega\times Y^*).
\endeq

\end{teo}

\begin{teo}\label{t2} Let $\di (u^{\ep})_{\ep}$ be a bounded 
sequence in $H^2(\Omega^\ep).$ Then there exists $u\in H^2(\Omega)$ and
$u^2\in L^2(\Omega;H^2_{\rm per}(Y)/{\RR})$ such that
\neweq{asapte}
\begin{tabular}{cll}
$\di {\cal T}^{\ep}(u^{\ep})\we u$& $\quad\mbox{weakly in }\quad$&$ L^2(\Omega\times Y^*),$\\
$\di {\cal T}^{\ep}(\nabla_x u^{\ep})\we \nabla_x u$&$\quad\mbox{weakly in }\quad$&
$L^2(\Omega\times Y^*),$\\
$\di {\cal T}^{\ep}(\nabla^2_xu^{\ep})\we \nabla^2_xu+\nabla^2_yu^2$&
$\quad\mbox{weakly in }\quad$& $L^2(\Omega\times Y^*).$\\
\end{tabular}
\endeq
\end{teo}

\noindent{\bf Proof.}
Since $(u^{\ep})$ is bounded in $H^2(\Omega^\ep),$ it follows that it weakly 
converges to some $u\in H^2(\Omega).$ According to Theorem \ref{t1}, there exists
$\di  u^1\in L^2(\Omega;H^1_{\rm per}(Y^*))$ such that
\neweq{aopt}
{\cal T}^{\ep}(u^{\ep})\we u\quad\mbox{ weakly in }\,L^2(\Omega\times Y^*),
\endeq
and
\neweq{anoua}
{\cal T}^{\ep}(\nabla_x u^{\ep})\we \nabla_x u+
\nabla_y u^1
\quad\mbox{ weakly in }\,L^2(\Omega\times Y^*).
\endeq
Moreover, by the boundedness of $(u^{\ep})$ in 
$H^2(\Omega^\ep),$
it follows that ${\cal T}^{\ep}(\nabla^2_xu^{\ep})$ is bounded in $L^2(\Omega\times Y^*).$
Hence, there exists $\varrho\in L^2(\Omega\times Y^*)$ such that
\neweq{azece}
{\cal T}^{\ep}(\nabla^2_x u^{\ep})\we \varrho
\quad\mbox{ weakly in }\,L^2(\Omega\times Y^*).
\endeq 
 Let $\psi\in{\cal D}(\Omega;C^{\infty}_{\rm per}(Y^*)).$ Then we have
\neweq{aunspe}
\di \int\limits_{\Omega}\frac{\partial^2u^{\ep}}{\partial x_i\partial x_j}\;
\psi\left(x,\frac{x}{\ep}\right)dx=
-\di \int\limits_{\Omega}\frac{\partial u^{\ep}}{\partial x_j}
\left[\frac{\partial \psi}{\partial x_i}+\frac{1}{\ep}\frac{\partial \psi}{\partial y_i}\right]
\left(x,\frac{x}{\ep}\right)dx
\endeq
Using the unfolding operator and \eq{atrei} in the above relation we have
\neweq{adoispe}
\di \ep\int\limits_{\Omega\times Y^*}{\cal T}^{\ep}\left(\frac{\partial^2u^{\ep}}
{\partial x_i\partial x_j}\right)\;\psi\left(x,y\right)dxdy=
-\di \int\limits_{\Omega\times Y^*}{\cal T}^{\ep}\left(
\frac{\partial u^{\ep}}{\partial x_j}\right)
\left[\ep\frac{\partial \psi}{\partial x_i}+\frac{\partial \psi}{\partial y_i}\right]
\left(x,y\right)dxdy.
\endeq
Passing to the limit in \eq{adoispe} and using 
\eq{anoua}, \eq{azece} we get
$$\di 0=\int\limits_{\Omega\times Y^*}\left\{
\frac{\partial u}{\partial x_j}+
\frac{\partial  u^1}{\partial y_j}\right\}\frac{\partial\psi}{\partial y_i}\left(x,y\right)dxdy
\quad\forall\,\psi\in{\cal D}(\Omega;C^{\infty}
_{\rm per}(Y^*)).$$
This yields that $\di 
\frac{\partial u}{\partial x_j}+
\frac{\partial u^1}{\partial y_j}$ does not depend on $y.$ 
Since $ u^1$ is $Y^*-$periodic in the second variable we conclude that 
$u^1(x,y)=u^1(x)$
and by \eq{anoua} we deduce  
$$\di {\cal T}^{\ep}(\nabla_x u^{\ep})\we \nabla_x u\quad\mbox{weakly in }\quad
L^2(\Omega\times Y^*).$$
\medskip

Let now $\di \psi\in{\cal D}(\Omega;C^{\infty}
_{\rm per}(Y^*))$ with $\di\nabla_y\psi(x,y)=0.$
From \eq{aunspe} we obtain
$$\di \int\limits_{\Omega}\frac{\partial^2u^{\ep}}{\partial x_i\partial x_j}\;
\psi\left(x,\frac{x}{\ep}\right)dx=
-\di \int\limits_{\Omega}\frac{\partial u^{\ep}}{\partial x_j}
\frac{\partial \psi}{\partial x_i}
\left(x,\frac{x}{\ep}\right)dx.$$
Using again the unfolding operator we have
\neweq{atreispe}
\di \int\limits_{\Omega\times Y^*}{\cal T}^{\ep}\left(\frac{\partial^2u^{\ep}}
{\partial x_i\partial x_j}\right)\;\psi\left(x,y\right)dxdy=
-\di \int\limits_{\Omega\times Y^*}{\cal T}^{\ep}
\left(\frac{\partial u^{\ep}}{\partial x_j}\right)\frac{\partial \psi}{\partial x_i}(x,y)dxdy.
\endeq
Passing to the limit we get
$$\begin{tabular}{lll}
$\di \int\limits_{\Omega\times Y^*}\varrho _{ij}(x,y)\psi(x,y)dxdy$&$
=$&$-\di \int\limits_{\Omega\times Y^*}
\frac{\partial u}{\partial x_j}(x,y)\frac{\partial\psi}{\partial x_i}
(x,y)dxdy$\\
&$=$&$\di \int\limits_{\Omega\times Y}
\frac{\partial^2 u}{\partial x_i\partial x_j}(x,y)\psi(x,y)dxdy.$\\
\end{tabular}$$
The above relations lead us to
$$\di \int\limits_{\Omega\times Y}\left[\varrho _{ij}(x,y)
-\frac{\partial^2 u}{\partial x_i\partial x_j}(x,y)\right]\psi(x,y)dxdy=0,$$
for all $\di \psi\in{\cal D}(\Omega;C^{\infty}
_{\rm per}(Y^*))$ with $\di\nabla_y\psi(x,y)=0.$\\
Then, there exists 
$\tilde u\in \left[L^2(\Omega;H^1_{\rm per}(Y^*)/\RR)\right]^N$ such that
\neweq{apaispe}
\di \varrho-\nabla_x^2u=\nabla_y\tilde u.
\endeq
Taking into account the symmetry of the left side member in the above equality, we
may conclude that there exists  
$u^2\in L^2(\Omega;H^2_{\rm per}(Y^*)/\RR)$ such that
$\di \tilde u=\nabla_yu^2.$ Now \eq{apaispe} becomes 
$$\di \varrho=\nabla^2_xu+\nabla^2_yu^2.$$
The proof of Theorem \ref{t2} is now complete.$\qed$
\medskip

We introduce now the averaging operator
$${\cal U}^\ep: L^2(\Omega^\ep\times Y^*)\ri
 L^2(\Omega^\ep),$$
$$ {\cal U}^\ep(\Phi)(x)=\frac{1}{|Y^*|}\int
_{Y^*}
\Phi\left(\ep\left[\frac{x}{\ep}\right]_{Y^*}+\ep z,
\left\{\frac{x}{\ep}\right\}_{Y^*}
\right)dz.$$

\begin{prop}\label{p3}{\rm (see \cite{Cior})}\\
{\rm (i)}\quad ${\cal U}^\ep(\phi)\ri\phi$ strongly in 
$L^2(\Omega),$ for all $\phi\in 
L^2(\Omega);$\\
{\rm (ii)}\quad $\di {\cal U}^\ep({\cal T}^\ep(\phi))=\phi,$ for all $\phi\in L^2(\Omega^\ep);$\\
{\rm (iii)}\quad 
$\di {\cal T}^\ep({\cal U}^\ep(\Phi))(x,y)
=\frac{1}{|Y^*|}\int_{Y^*}
\Phi\left(\ep\left[\frac{x}{\ep}\right]_{Y^*}+\ep z,
\left\{\frac{x}{\ep}\right\}_{Y^*}
\right)dz,$ for all $\Phi\in
L^2(\Omega^\ep\times Y^*);$\\
{\rm (iv)}\quad
$\di \int_{\Omega^\ep}{\cal U}^\ep(\Phi)(x)dx
=\frac{1}{|Y^*|}\int
_{\Omega^\ep\times Y^*}
\Phi(x,y)dxdy,$ for all 
$\Phi\in
L^2(\Omega^\ep\times Y^*);$\\
{\rm (v)}\quad
$\di {\cal U}^\ep(\Phi)\ri \frac{1}{|Y^*|}\int
_{Y^*}\Phi(x,y)dy,$ for all 
$\Phi\in L^2(\Omega\times Y^*).$\\
\end{prop}
\begin{teo}\label{topu}{\rm (see \cite{Cior})}
Let $\di (\phi^{\ep})_{\ep}\subset L^2(\Omega).$
The following weak convergences are equivalent\\
{\rm (i)}\quad ${\cal T}^\ep(\phi^\ep)\we\Phi$ in $L^2(\Omega\times Y^*);$\\
{\rm (ii)}\quad 
$\phi^\ep-{\cal U}^\ep(\Phi)\we\phi$ in 
$L^2(\Omega).$\\
A similar equivalence holds for strong convergences.
\end{teo}
\section{The membrane problem}
\subsection{The convergence results}
Latin indices take their values in the set $\{1,2,3\}$ and Greek indices take their 
values in $\{1,2\}.$ The sumation convention is also used. Boldface letters 
represent vector-valued functions.
\smallskip

Let $\omega\subset\RR^2$ be an open bounded and connected set with a 
Lipschitz-continuous boundary. Let $\varphi:\overline{\omega}\ri\RR^3$ be a 
$C^2(\overline{\omega})$ one-to-one function such that the vectors 
$\di {\bf a_1}=\partial\varphi/\partial x_1$ and 
$\di {\bf a_2}=\partial\varphi/\partial x_2$ are linearly independent at all points
in $\overline{\omega}.$ Denote 
$$\di {\bf a_3}=\frac{{\bf a_1}\times{\bf a_2}}{|{\bf a_1}\times{\bf a_2}|}
\quad\mbox{and }\quad a=\det({\bf a}_{\alpha}\cdot {\bf a}_{\beta}).$$

We consider the membrane problem
\neweq{doiunu}
\left\{\begin{tabular}{ll}
$\di \int_{\omega^{\ep}}c^{\ep}({\bf u}^{\ep},{\bf  v})+
e^{\ep}(\varphi^{\ep},{\bf v})$&$\di =
\int_{\omega^{\ep}}F_iv_i\sqrt{a^{\ep}}+
\int_{\Gamma_+^{\ep}\cup\Gamma_-^{\ep}}q^iv_i
\sqrt{a^{\ep}}+
\int_{\omega^{\ep}} h^{\alpha\beta,\ep}\gamma_{\alpha\beta}({\bf v})\sqrt{a^{\ep}}
$\\
$\di\int_{\omega^{\ep}}e^{\ep}(\psi,{\bf u}^{\ep})+d^{\ep}(\varphi^{\ep},\psi)$
&$\di =
\int_{\omega^{\ep}}h^{\alpha}\partial_\alpha\psi\sqrt{a^{\ep}},$\\
\end{tabular}\right.\endeq
for all $({\bf v},\psi)\in {\bf V}(\omega^\ep)\times W(\omega^\ep)$,
where
$$\begin{tabular}{ll}
$\di {\bf V}(\omega^\ep)$&$\di=\{{\bf v}\in H^1(\omega^\ep)\times
 H^1(\omega^\ep)\times L^2(\omega^\ep);\;\;v_\alpha=0\;\;\mbox{ on }\;
\gamma_0^M \},\;\;\;\;\mbox{ with } meas\; \gamma_0^M>0,$\\
$\di W(\omega^\ep)$&
$\di=\{\psi\in H^1(\omega^\ep)
;\;\;\psi=0\;\;\mbox{ in }\omega^\ep\setminus\omega_0,
\psi=0\;\;\mbox{ on }\;\gamma_0^E \}.$\\
\end{tabular}$$
$$\begin{tabular}{ll}
$c^{\ep}({\bf u},{\bf v})$&
$=c^{\alpha\beta\lambda\mu,\ep}\sqrt{a^{\ep}}\;\gamma_{\alpha\beta}({\bf u})
\gamma_{\lambda\mu}({\bf v}),$\\
$e^{\ep}({\bf v},\psi)$&
$=e^{\lambda\alpha\beta,\ep}\sqrt{a^{\ep}}\;\partial_\lambda\psi
\gamma_{\alpha\beta}({\bf v}),$\\
$d^{\ep}(\psi,\varphi)$&
$=d^{\alpha\lambda,\ep}\sqrt{a^{\ep}}\;\partial_\alpha\psi
\partial_\lambda\varphi,$\\
$\di \gamma_{\alpha\beta}({\bf v})$
&$=s_{\alpha\beta}({\bf v})-\Gamma^k_{\alpha\beta}v_k,$\\
$ s_{\alpha\beta}({\bf v})$&
$=\frac{1}{2}\left(\partial_\alpha v_\beta+
\partial_\beta v_\alpha\right),$\\
\end{tabular}$$
$$F^i=F^i(x_1,x_2)=\int_{-1}^1f^{i}(x_1,x_2,t)dt,\quad
\mbox{ where }\;\;{\bf f}=(f^i)\in {\bf L}^2(\omega\times [-1,1]),$$
$$ {\bf q}=(q^i)\in{\bf L}^2(\Gamma_{+}\cup\Gamma_{-}).$$

Taking into account the ellipticity and the symmetry of the tensors 
$(c^{\alpha\beta\lambda\mu})$
and $(d^{\alpha\lambda})$ we can apply the
Lax-Milgram Theorem in order to get the existence of a 
unique displacement ${\bf u}^{\ep}$ and of a unique electric potential
$\varphi^{\ep}$ that
 verifies the variational equation \eq{doiunu}.
\smallskip

\begin{teo}\label{t3}
The sequences  $\left({\cal T}^{\ep}({\bf u}^{\ep})\right)_{\ep}$ and
$\left({\cal T}^{\ep}(\varphi^{\ep})\right)_{\ep}$
weakly converge to
${\bf u}\in{\bf V}(\omega)$ and  $\varphi\in W(\omega)$ respectively,
which are the unique solutions of the
homogenized problem 
\neweq{doiec}
\left\{\begin{tabular}{ll}
$\di\int_{\omega\times Y^*}\left[{\overline c}({\bf u},{\bf v})+
{\overline e}({\bf v},\varphi)\right]dxdy$&$\di =
|Y^*|_a\int_{\omega}F^iv_{i}dx+|Y^*|_a\int_{\Gamma_{-}\cup\Gamma_{+}}
q^iv_{i}d\Gamma, $\\
$\di\int_{\omega\times Y^*}[-{\overline e}({\bf u},{\bf v})+
{\overline d}(\psi,\varphi)]dxdy$&$=0,$\\
\end{tabular}\right.
\endeq
for all $\;({\bf v},\psi)\in{\bf V}(\omega)\times W(\omega),$ where
$\di |Y^*|_a=\int_{Y^*}\sqrt{a}dy$ and
$$\left\{\begin{tabular}{ll}
${\overline c}({\bf u},{\bf v})$&$=
{\overline c}^{\alpha\beta\lambda\mu}\sqrt{a}\;\gamma_{\alpha\beta}({\bf u})\;
\gamma_{\lambda\mu}({\bf v}),$\\
${\overline e}({\bf v},\psi)$&
$={\overline e}^{\lambda\alpha\beta}
\sqrt{a}\;\partial_\lambda\psi\;\gamma_{\alpha\beta}({\bf v}),$\\
$\di {\overline d}(\psi,\varphi)$&
$\di ={\overline d}^{\alpha\lambda}\sqrt{a}\;
\partial_\alpha\psi\;\partial_\lambda\varphi,$\\
\end{tabular}\right.$$
for all ${\bf u},{\bf v}\in{\bf H}^1_{\rm per}(Y^*)/\RR,$
$\varphi,\psi\in H^1_{\rm per}(Y^*)/\RR,$
\neweq{xcoef}
\left\{\begin{tabular}{ll}
$\di {\overline c}^{\alpha\beta\tau\theta}$&$\di =
\int_{Y^*}\left[c^{\alpha\beta\lambda\mu}\sqrt{a}\,
(\delta^{\tau\theta}_{\lambda\mu}+s_{\lambda\mu,y}({\bf w}^{\tau\theta}
))+e^{\lambda\alpha\beta}\sqrt{a}\partial_{\lambda,y}
\zeta^{\tau\theta}\right]dy,$\\
$\di {\overline e}^{\sigma\alpha\beta}$&$\di=
\int_{Y^*}\left[c^{\alpha\beta\lambda\mu}\sqrt{a}\,s_{\lambda\mu}(z^\sigma)+
e^{\lambda\alpha\beta}\sqrt{a}
(\delta^{\sigma}_{\lambda}+\partial_{\lambda,y}\eta^\sigma)
\right]dy,$\\
$ \di {\overline d}^{\alpha\sigma}$&$\di =
\int_{Y^*}\left[-e^{\alpha\lambda\mu}\sqrt{a}\,
s_{\lambda\mu,y}(z^\sigma)+d^{\alpha\lambda}\sqrt{a}(\delta^\sigma_\lambda
+\partial_{\lambda,y}\eta^{\sigma})\right]dy.$\\
\end{tabular}\right.\endeq
The local functions $({\bf w}^{\tau\theta},\zeta^{\tau\theta})$
and $({\bf z}^{\sigma},\eta^{\sigma})$
verifies the local problems 
\neweq{local1}
\left\{\begin{tabular}{ll}
$\di
\int_{Y^*}\left(c_y(\Sigma^{\tau\theta}+{\bf w}^{\tau\theta},{\bf v})+
e_y({\bf v},\zeta^{\tau\theta})
\right)dy$&$=0$,\\
$\di\int_{Y^*}\left(-e_y(\Sigma^{\tau\theta}+{\bf w}^{\tau\theta},\psi)
+d_y(\zeta^{\tau\theta},\psi)\right)dy
$&$=0,$\\
\end{tabular}\right.\endeq
\neweq{local2}
\left\{\begin{tabular}{ll}
$\di \int_{Y^*}\left(c_y({\bf z}^\sigma,{\bf v})+e_y({\bf v},y_\sigma+\eta^\sigma)
\right)dy$&$=0,$\\
$\di \int_{Y^*}\left(-e_y({\bf z}^\sigma,\psi)+d_y(y_\sigma+\eta^{\sigma},
\psi)\right)dy$&$=0,$\\
\end{tabular}\right.\endeq
for all $({\bf v},\psi)\in  {\bf H}^1_{\rm per}(Y^*)/\RR \times 
H^1_{\rm per}(Y^*)/\RR,$
where
$\Sigma^{\alpha\beta}=y_\alpha e_\beta+y_\beta e_\alpha$ and 
\begin{eqnarray}\label{locy1}
\di c_y({\bf w},{\bf v})&=&c^{\alpha\beta\tau\theta}\sqrt{a}\,
s_{\alpha\beta,y}({\bf w}) s_{\lambda\mu,y}({\bf v})
,\\
\label{locy2}
\di e_y({\bf v},\psi)&=& e^{\lambda\alpha\beta}\sqrt{a}\;
 s_{\alpha\beta}({\bf v})\partial_{\lambda,y}(\psi),\\
\label{locy3}
\di d_y(\psi,\phi)&=& d^{\alpha\lambda}\sqrt{a}\,
\partial_{\alpha,y}\psi\;\partial_{\lambda,y}\phi.
\end{eqnarray}
\end{teo}
\smallskip

\noindent{ \bf Proof.}
To obtain an {\it a priori estimate}, we take
 ${\bf v}={\bf u}^{\ep}$ and $\psi=\varphi^{\ep}$ in \eq{doiunu}. 
By Korn and Poincar\'e's inequalities for perforated domains, it follows
$$\|{\bf u}^{\ep}\|_{{\bf V}(\omega^\ep)}+\|\varphi^{\ep}\|_{W(\omega^\ep)}\leq C,$$
where $C$ is a positive constant that depends only on $\omega$ (but not on $\ep$).

By Theorem \ref{t1}, there exist $({\bf u},\varphi)\in{\bf V}(\omega)\times W(\omega)$ and
two fields of corectors 
${\bf u}^1=(u_1^1,u_2^1)\in {\bf  L}^2(\omega;H^1_{\rm per}(Y)/\RR),$
 $\varphi^1\in L^2(\omega;H^1_{\rm per}(Y^*))$ such that, up to a 
sequence we have
\neweq{doisapte}
{\cal T}^{\ep}({\bf u}^{\ep})\we {\bf u}\quad\mbox{ weakly in }\,
{\bf L}^2(\omega\times Y^*;\RR^3),
\endeq
\neweq{doiopt}
{\cal T}^{\ep}(\nabla_x ({\bf u}^{\ep})\we \nabla_x({\bf u})+\nabla_y {\bf u}^1
\quad\mbox{ weakly in }\,{\bf L}^2(\omega\times Y^*;\RR^2),
\endeq
\neweq{doioptbis}
{\cal T}^{\ep}(\nabla_x ({\varphi}^{\ep})\we \nabla_x({\varphi})+\nabla_y {\varphi}^1
\quad\mbox{ weakly in }\,{\bf L}^2(\omega\times Y^*;\RR^2).
\endeq
The linerity of ${\cal T}^{\ep}$ implies
\neweq{doinoua}
{\cal T}^{\ep}(\gamma_{\alpha\beta}({\bf u}^{\ep}))\we \gamma_{\alpha\beta,x}({\bf u})+
s_{\alpha\beta,y}({\bf u}^1)
\quad\mbox{ weakly in }\,L^2(\omega\times Y^*;\RR^2)
\endeq
and
\neweq{doizece}
{\cal T}^{\ep}(\partial_{\alpha}\varphi^{\ep})\we 
\partial_{\alpha,x}\varphi+\partial_{\alpha,y}\varphi^1
\quad\mbox{ weakly in }\,L^2(\omega\times Y^*;\RR^2).
\endeq

We chose as a test function in \eq{doiunu}
$$\di {\bf v}^\ep(x)={\bf v}_1(x)+\ep{\bf v}_2\left(x,\frac{x}{\ep}\right),
\quad {\bf v}_1\in{\cal D}(\omega);{\bf v}_2\in{\cal D}(\omega;
C^{\infty}_{\rm per}(Y^*)),$$

$$\di \psi^\ep(x)=\psi_1(x)+\ep\psi_2\left(x,\frac{x}{\ep}\right),
\quad \psi_1\in{\cal D}(\omega);\psi_2\in{\cal D}(\omega;C^{\infty}_{\rm per}(Y^*)).$$
It follows that
\neweq{doidoispe}
\di {\bf v}^\ep(x)\ri{\bf v}_1(x)\quad\mbox{ strongly in }\,{\bf L}^2(\omega),\endeq
$$\di \nabla_x{\bf v}^\ep(x)\we\nabla_x{\bf v}_1(x)+\nabla_y{\bf v}_2\left(x,y\right)
\quad\mbox{ weakly in }\,{\bf L}^2(\omega\times Y^*),$$
\neweq{doitreispe}
{\cal T}^{\ep}(\gamma_{\rho\sigma}({\bf v}^\ep))\we 
\gamma_{\rho\sigma}(\phi)+e_{\rho\sigma,y}(\psi)
\quad\mbox{ weakly in }\,L^2(\omega\times Y^*),
\endeq
\neweq{doitreispebis}
{\cal T}^{\ep}(\partial_\alpha\psi^\ep)\we\partial_{\alpha,x}\psi_1 
+\partial_{\alpha,y}\psi^1\quad\mbox{ weakly in }\,L^2(\omega\times Y^*).
\endeq

Using the linearity of ${\cal T}^{\ep}$ and the fact that
$\di \int_{\omega}v=\frac{1}{|Y^*|}\int_{\omega\times Y^*}{\cal T}^{\ep}(v)$
for all $v\in L^2(\omega),$ we can pass to the limit in the variational form
\eq{doiunu}. We get
\neweq{doipaispe}
\left\{\begin{tabular}{ll}
$\di\int_{\omega\times Y^*}\left[c^{\alpha\beta\gamma\lambda}\sqrt{a}\left(
\gamma_{\lambda\mu,x}({\bf u})+s_{\lambda\mu,y}({\bf u}^1\right)+
e^{\lambda\alpha\beta}\sqrt{a}
\left(\partial_\lambda\varphi+\partial_{\lambda,y}\varphi^1\right)\right]
\left(\gamma_{\alpha\beta,x}({\bf v}_1)+s_{\alpha\beta,y}({\bf v}_2\right)
dxdy$\\
$\di\quad=|Y^*|_a\int_{\omega}F^iv_{1i}dx+|Y^*|_a\int_{\Gamma_{-}\cup\Gamma_{+}}
q^iv_{1i}d\Gamma+
\int_{\omega\times Y^*}h^{\alpha\beta}\sqrt{a}\left[
(\gamma_{\alpha\beta,x}({\bf v}_1)+s_{\alpha\beta,y}({\bf v}_2)\right]dxdy,$\\
$\di\int_{\omega\times Y^*}\left[-e^{\alpha\lambda\mu}\sqrt{a}\left(
\gamma_{\lambda\mu,x}({\bf u})+s_{\lambda\mu,y}({\bf u}^1\right)+
d^{\alpha\lambda}\sqrt{a}
\left(\partial_\lambda\varphi+\partial_{\lambda,y}\varphi^1\right)\right]
\left(\partial_{\alpha,x}\psi_1+\partial_{\alpha,y}\psi_2\right)
dxdy$\\
$ \di\quad =
\int_{\omega\times Y^*}l^{\alpha}\sqrt{a}
\left(\partial_{\alpha,x}\psi_1+\partial_{\alpha,y}\psi_2\right)dxdy$\\
\end{tabular}\right.\endeq

Letting ${\bf v}_2=0,$ $\psi_2=0$ in the above equalities, by density
 we have
\neweq{doix1}
\left\{\begin{tabular}{ll}
$\di\int_{\omega\times Y^*}\left[c^{\alpha\beta\gamma\lambda}\sqrt{a}\left(
\gamma_{\lambda\mu,x}({\bf u})+s_{\lambda\mu,y}({\bf u}^1\right)+
e^{\lambda\alpha\beta}\sqrt{a}
\left(\partial_\lambda\varphi+\partial_{\lambda,y}\phi^1\right)\right]
\gamma_{\alpha\beta,x}({\bf v})dxdy$\\
$\di=|Y^*|_a\int_{\omega}F^iv_{i}dx+|Y^*|_a\int_{\Gamma_{-}\cup\Gamma_{+}}
q^iv_{i}d\Gamma+
\int_{\omega\times Y^*}h^{\alpha\beta}\sqrt{a}
\gamma_{\alpha\beta,x}({\bf v})dxdy$\\
$\di\int_{\omega\times Y^*}\left[-e^{\alpha\lambda\mu}\sqrt{a}\left(
\gamma_{\lambda\mu,x}({\bf u})+s_{\lambda\mu,y}({\bf u}^1\right)+
d^{\alpha\lambda}\sqrt{a}
\left(\partial_\lambda\varphi+\partial_{\lambda,y}\phi^1\right)\right]
\partial_{\alpha,x}\psi dxdy$\\
$\di=
\int_{\omega\times Y^*}l^{\alpha}\sqrt{a}
\partial_{\alpha,x}\psi dxdy,$\\
\end{tabular}\right.\endeq
for all $({\bf v},\psi)\in{\bf V}(\omega)\times W(\omega).$
By linearity, we take
\neweq{usapo}
\di {\bf u}^1(x,y)=\gamma_{\tau\theta}({\bf u}(x,y)){\bf w}^{\tau\theta}(y)+
\partial_\sigma\varphi(x){\bf z}^\sigma(y)+q(x,y),
\endeq
\neweq{fisapo}
\varphi^1(x,y)=
\gamma_{\tau\theta}({\bf u}(x,y))\zeta^{\tau\theta}(y)+
\partial_\sigma\varphi(x)\eta^\sigma(y)+\xi(x,y),
\endeq
where ${\bf w}^{\tau\theta},\zeta^{\tau\theta},{\bf z}^\sigma,\eta^\sigma$
are $Y^*$ periodic in $y$ and
\neweq{doisaptispe}
c^{\alpha\beta\gamma\lambda}s_{\lambda\mu,y}(q(x,y))+
e^{\lambda\alpha\beta}\partial_{\lambda,y}\xi(x,y)=h^{\alpha\beta}(x,y),
\endeq
\neweq{doioptispe}
-e^{\alpha\lambda\mu}s_{\lambda\mu,y}(q(x,y))+
d^{\alpha\lambda}\partial_{\lambda,y}\xi(x,y)=l^\alpha(x,y).
\endeq
Replacing \eq{usapo} and \eq{fisapo} in \eq{doix1} and
 taking into account the properties of $q(x,y)$ and $\xi(x,y)$ we obtain
\neweq{doinouaspe}
\left\{\begin{tabular}{ll}
$\di\int_{\omega\times Y^*}\left[{\overline c}^{\alpha\beta\tau\theta}
\gamma_{\tau\theta}({\bf u})+
{\overline e}^{\sigma\alpha\beta}\partial_\sigma\varphi\right]
\gamma_{\alpha\beta}({\bf v})dxdy=
|Y^*|_a\int_{\omega}F^iv_{i}dx+|Y^*|_a\int_{\Gamma_{-}\cup\Gamma_{+}}
q^iv_{i}d\Gamma $\\
$\di\int_{\omega\times Y^*}[-{\overline  f}^{\alpha\tau\theta}
\gamma_{\tau\theta}({\bf u})+{\overline d}^{\alpha\sigma}
\partial_\lambda\varphi]\partial_{\alpha}\psi dxdy=0,$\\
\end{tabular}\right.
\endeq
for all $({\bf v},\psi)\in{\bf V}(\omega)\times W(\omega),$
where
\begin{eqnarray}\label{doidouazeci}
\di {\overline c}^{\alpha\beta\tau\theta}&=&
\int_{Y^*}\left[c^{\alpha\beta\lambda\mu}\sqrt{a}\,
(\delta^{\tau\theta}_{\lambda\mu}+s_{\lambda\mu,y}({\bf w}^{\tau\theta}
))+e^{\lambda\alpha\beta}\sqrt{a}\partial_{\lambda,y}\zeta^{\tau\theta}\right]dy,\\
\label{doidouazeciunu}
\di {\overline e}^{\sigma\alpha\beta}&=&
\int_{Y^*}\left[c^{\alpha\beta\lambda\mu}\sqrt{a}\,s_{\lambda\mu}
({\bf z}^\sigma)+
e^{\lambda\alpha\beta}\sqrt{a}
(\delta^{\sigma}_{\lambda}+\partial_{\lambda,y}\eta^\sigma)
\right]dy,\\
\label{doidouazecidoi}
\di {\overline f}^{\alpha\tau\theta}&=&
\int_{Y^*}\left[e^{\alpha\lambda\mu}\sqrt{a}\,
(\delta^{\tau\theta}_{\lambda\mu}+s_{\lambda\mu,y}({\bf w}^{\tau\theta}))
-d^{\alpha\lambda}\sqrt{a}\partial_{\lambda,y}\zeta^{\tau\theta}\right]dy,\\
\label{doidouazecitrei}
\di {\overline d}^{\alpha\sigma}&=&
\int_{Y^*}\left[-e^{\alpha\lambda\mu}\sqrt{a}\,
s_{\lambda\mu,y}({\bf z}^\sigma)+d^{\alpha\lambda}\sqrt{a}(\delta^\sigma_\lambda
+\partial_{\lambda,y}\eta^{\sigma})\right]dy.
\end{eqnarray}
We will prove that
${\overline e}^{\sigma\alpha\beta}={\overline f}^{\sigma\alpha\beta}.$ In
that sense we first determine the local problems verified by the functions
 ${\bf w}^{\tau\theta},\zeta^{\tau\theta},{\bf z}^\sigma,\eta^\sigma.$
From \eq{doidouazeci}-\eq{doidouazecitrei} we have
\begin{eqnarray}\label{doidouazecipatru}
\di {\overline c}^{\alpha\beta\tau\theta}&=&
\int_{Y^*}\left[c^{\alpha\beta\lambda\mu}\sqrt{a}\,
s_{\lambda\mu,y}(\Sigma^{\tau\theta}+{\bf w}^{\tau\theta})+
e^{\lambda\alpha\beta}\sqrt{a}
\partial_{\lambda,y}\zeta^{\tau\theta}\right]dy,\\
\label{doidouazecicinci}
\di {\overline e}^{\sigma\alpha\beta}&=&
\int_{Y^*}\left[c^{\alpha\beta\lambda\mu}\sqrt{a}\,s_{\lambda\mu}(z^\sigma)+
e^{\lambda\alpha\beta}\sqrt{a}\;
\partial_{\lambda,y}(y_{\sigma}+\eta^\sigma)
\right]dy,\\
\label{doidouazesase}
\di {\overline f}^{\alpha\tau\theta}&=&
\int_{Y^*}\left[e^{\alpha\lambda\mu}\sqrt{a}\,
s_{\lambda\mu,y}((\Sigma^{\tau\theta}+{\bf w}^{\tau\theta}))
-d^{\alpha\lambda}\sqrt{a}\partial_{\lambda,y}\zeta^{\tau\theta}\right]dy,\\
\label{doidouazecisapte}
\di {\overline d}^{\alpha\sigma}&=&
\int_{Y^*}\left[-e^{\alpha\lambda\mu}\sqrt{a}\,
s_{\lambda\mu,y}(z^\sigma)+d^{\alpha\lambda}\sqrt{a}\;\partial_{\lambda,y}
(y_\sigma+\eta^{\sigma})\right]dy,
\end{eqnarray}

We let ${\bf v}_1=0$ and $\psi_1=0$ in \eq{doipaispe}.
 By density we deduce
\neweq{doidouazecipbb}
\left\{\begin{tabular}{ll}
$\di\int_{\omega\times Y^*}\left[c^{\alpha\beta\gamma\lambda}\sqrt{a}\left(
\gamma_{\lambda\mu,x}({\bf u})+s_{\lambda\mu,y}({\bf u}^1)\right)+
e^{\lambda\alpha\beta}\sqrt{a}\;
\left(\partial_\lambda\varphi+\partial_{\lambda,y}\varphi^1\right)\right]
s_{\alpha\beta,y}({\bf v})dxdy$\\
$\qquad\di=\int_{\omega\times Y^*}h^{\alpha\beta}\sqrt{a}\;
s_{\alpha\beta,y}({\bf v})dxdy,$\\
$\di\int_{\omega\times Y^*}\left[-e^{\alpha\lambda\mu}\sqrt{a}\left(
\gamma_{\lambda\mu,x}({\bf u})+s_{\lambda\mu,y}({\bf u}^1)\right)+
d^{\alpha\lambda}\sqrt{a}
\left(\partial_\lambda\varphi+\partial_{\lambda,y}\varphi^1\right)\right]
\partial_{\alpha,y}\psi dxdy$\\
$\qquad\di=\int_{\omega\times Y^*}l^{\alpha}\sqrt{a}
\partial_{\alpha,y}\psi dxdy,$\\
\end{tabular}\right.\endeq
for all $({\bf v},\psi)\in {\bf L}^2(\omega;{\bf H}^1_{\rm per}(Y^*)/ \RR)
\times  L^2(\omega; H^1_{\rm per}(Y^*)/ \RR).$\\
Using again \eq{usapo} and \eq{fisapo} we deduce
\neweq{doidouazecisase}
\left\{\begin{tabular}{ll}
$\di
\int_{Y^*}\left(c_y(\Sigma^{\tau\theta}+{\bf w}^{\tau\theta},{\bf v})+
e_y({\bf v},\zeta^{\tau\theta})
\right)dy$&$=0$,\\
$\di\int_{Y^*}\left(-e_y(\Sigma^{\tau\theta}+{\bf w}^{\tau\theta},\psi)
+d_y(\zeta^{\tau\theta},\psi)\right)dy
$&$=0,$\\
\end{tabular}\right.\endeq

\neweq{doidouazecisapte}
\left\{\begin{tabular}{ll}
$\di \int_{Y^*}\left(c_y({\bf z}^\sigma,{\bf v})+e_y({\bf v},y_\sigma+\eta^\sigma)
\right)dy$&$=0,$\\
$\di \int_{Y^*}\left(-e_y({\bf z}^\sigma,\psi)+d_y(y_\sigma+\eta^{\sigma},
\psi)\right)dy$&$=0,$\\
\end{tabular}\right.\endeq
for all $({\bf v},\psi)\in  {\bf H}^1_{\rm per}(Y^*)/\RR \times 
H^1_{\rm per}(Y^*)/\RR,$
where $c_y,e_y,d_y$ are given by \eq{locy1}-\eq{locy3}.

Next, we prove that ${\overline e}^{\sigma\tau\theta}=
{\overline f}^{\sigma\tau\theta}.$
By taking $({\bf v},\psi)=({\bf z}^\sigma,\eta^\sigma)$ in 
\eq{doidouazecisase}, we get
$$\left\{\begin{tabular}{ll}
$\di
\int_{Y^*}\left(c_y({\bf w}^{\tau\theta},{\bf z}^\sigma)+
e_y({\bf z}^\sigma,\zeta^{\tau\theta})
\right)dy$&$\di =-\int_{Y^*}c^{\tau\theta\lambda\mu}\sqrt{a}\;
s_{\lambda\mu}({\bf z}^\sigma),$\\
$\di\int_{Y^*}\left(-e_y({\bf w}^{\tau\theta},\eta^\sigma)
+d_y(\zeta^{\tau\theta},\eta\sigma)\right)dy
$&$\di =\int_{Y^*}e^{\lambda\tau\theta}\sqrt{a}\partial_{\lambda,y}\eta^\sigma dy,$\\
\end{tabular}\right.$$
and so, by \eq{doidouazeciunu} it follows
\neweq{ee}
\di {\overline e}^{\sigma\tau\theta}=\int_{Y^*}\left(e^{\sigma\tau\theta}-
c_y({\bf w}^{\tau\theta},{\bf z}^\sigma))-e_y({\bf w}^{\tau\theta},\eta^\sigma)
-e_y({\bf z}^\sigma,\zeta^{\tau\theta})+
d_y(\eta^\sigma,\zeta^{\tau\theta})\right)dy.
\endeq
In the same manner, by taking $({\bf v},\psi)=({\bf z}^\sigma,\eta^\sigma)$ as a test
function in \eq{doidouazecisapte}, we get
$$\left\{\begin{tabular}{ll}
$\di
\int_{Y^*}\left(c_y({\bf w}^{\tau\theta},{\bf z}^\sigma)+
e_y({\bf w}^{\tau\theta},\eta^\sigma)
\right)dy$&$=\di -\int_{Y^*}e^{\sigma\lambda\mu}\sqrt{a}\;
s_{\lambda\mu}({\bf w}^{\tau\theta}),$\\
$\di\int_{Y^*}\left(-e_y({\bf z}^\sigma,\zeta^{\tau\theta})
+d_y(\eta^\sigma,\zeta^{\tau\theta})\right)dy
$&$=\di -\int_{Y^*}d^{\sigma\lambda}\sqrt{a}\;\partial_{\lambda,y}
\zeta^{\tau\theta} dy,
$\\
\end{tabular}\right.$$
and by \eq{doidouazecidoi} it follows
\neweq{ff}
\di {\overline f}^{\sigma\tau\theta}=\int_{Y^*}\left(e^{\sigma\tau\theta}-
c_y({\bf w}^{\tau\theta},{\bf z}^\sigma))-e_y({\bf w}^{\tau\theta},\eta^\sigma)
-e_y({\bf z}^\sigma,\zeta^{\tau\theta})+
d_y(\eta^\sigma,\zeta^{\tau\theta})\right)dy.\endeq
Now, by \eq{ee} and \eq{ff} we can conclude that
 ${\overline e}^{\sigma\tau\theta}={\overline f}^{\sigma\tau\theta}.$
\medskip

To prove that the homogenized problem is well posed, it remains only to show the symmetry and coercivity of $({\overline c}^{\alpha\beta\tau\theta})$ and
$({\overline d}^{\alpha\sigma}).$
\medskip

\noindent (i) \underline {The ellipticity and symmetry of the tensor}
 $({\overline c}^{\alpha\beta\tau\theta}).$
\smallskip

Let us prove first the symmetry. From \eq{doidouazeci} we deduce
$${\overline c}^{\alpha\beta\tau\theta}=
{\overline c}^{\beta\alpha\tau\theta}=
{\overline c}^{\alpha\beta\theta\tau}.$$
It remains to show that
${\overline c}^{\alpha\beta\tau\theta}=
{\overline c}^{\tau\theta\alpha\beta}.$
From \eq{doidouazecipatru} we have
\neweq{c0}
\di {\overline c}^{\alpha\beta\tau\theta}=
\int_{Y^*}\left[c_y(\Sigma^{\tau\theta}+{\bf w}^{\tau\theta},\Sigma^{\alpha\beta})
+e_y(\Sigma^{\alpha\beta},\zeta^{\tau\theta})\right]dy.
\endeq
But
\neweq{c1}
\di \int_{Y^*}e_y(\Sigma^{\alpha\beta},\zeta^{\tau\theta})dy=
-\int_{Y^*}e_y({\bf w}^{\alpha\beta},\zeta^{\tau\theta})dy+\int_{Y^*}
e_y(\Sigma^{\alpha\beta}+{\bf w}^{\alpha\beta},\zeta^{\tau\theta})
dy. \endeq
On the other hand,
by taking $({\bf v},\psi)=({\bf w}^{\alpha\beta},\zeta^{\alpha\beta})$ in 
\eq{doidouazecisase}, we get
\neweq{c2}
\left\{\begin{tabular}{ll}
$\di -\int_{Y^*}e_y({\bf w}^{\alpha\beta},\zeta^{\tau\theta})dy
$&$\di =\int_{Y^*}c_y(\Sigma^{\tau\theta}+{\bf w}^{\tau\theta},{\bf w}^{\alpha\beta}
)dy,$\\
$\di \int_{Y^*}e_y(\Sigma^{\alpha\beta}+{\bf w}^{\alpha\beta},\zeta^{\tau\theta})dy
$&$\di =\int_{Y^*}d_y(\zeta^{\alpha\beta},\zeta^{\tau\theta})dy.$\\
\end{tabular}\right.\endeq
Replacing now \eq{c1}-\eq{c2} in \eq{c0}
we deduce
$$\di {\overline c}^{\alpha\beta\tau\theta} =\int_{Y^*}
c_y( 
\Sigma^{\tau\theta}+{\bf w}^{\tau\theta},
\Sigma^{\alpha\beta}+{\bf w}^{\alpha\beta})dy+
\int_{Y^*}d_y(\zeta^{\alpha\beta},\zeta^{\tau\theta})dy.$$
From the above relations we can easily conclude the symmetry of the tensor
$({\overline c}^{\alpha\beta\tau\theta}).$
\medskip

Let us prove the coercivity.
Let $(X_{\alpha\beta})$ be a symmetric tensor (i.e. $X_{\alpha\beta}=
X_{\beta\alpha}$). First we note that by \eq{doidouazeci} we have
\neweq{doitreizecipatru}
\di {\overline c}^{\alpha\beta\tau\theta}X_{\alpha\beta}X_{\tau\theta}=
\int_{Y^*}c^{\alpha\beta\lambda\mu}\sqrt{a}(s_{\lambda\mu,y}({\bf W})
+X_{\lambda\mu})X_{\alpha\beta}+
\int_{Y^*}e^{\lambda\alpha\beta}\sqrt{a}(\partial_{\lambda,y}\Lambda)
X_{\alpha\beta}dy,
\endeq
where ${\bf W}={\bf w}^{\tau\theta}X_{\tau\theta}$ and 
$\Lambda=\zeta^{\tau\theta}X_{\tau\theta}.$
On the other hand, $({\bf W},\Lambda)$ is a solution of the following problem
\neweq{c4}
\left\{\begin{tabular}{ll}
$\di
\int_{Y^*}\left(c_y(X_{\tau\theta}\Sigma^{\tau\theta}+{\bf W},{\bf v})+
e_y({\bf v},\Lambda)
\right)dy$&$=0$,\\
$\di\int_{Y^*}\left(-e_y(X_{\tau\theta}\Sigma^{\tau\theta}+{\bf W},\psi)
+d_y(\Lambda,\psi)\right)dy
$&$=0,$\\
\end{tabular}\right.\endeq
for all $({\bf v},\psi)\in  {\bf H}^1_{\rm per}(Y^*)/\RR \times 
H^1_{\rm per}(Y^*)/\RR.$

Thus $({\bf W},\Lambda)$ is a saddle point of the following functional
$$\di I:{\bf H}^1_{\rm per}(Y^*)/\RR \times 
 H^1_{\rm per}(Y^*)/\RR\ri\RR,$$
defined by
$$\begin{tabular}{ll}
$\di I({\bf v},\psi)$&
$\di =\frac{1}{2}\int_{Y^*}
c^{\alpha\beta\lambda\mu}\sqrt{a}(s_{\lambda\mu,y}({\bf v})
+X_{\lambda\mu})(s_{\alpha\beta,y}({\bf v})+X_{\alpha\beta})dy$\\
&$\di \;\;\;+\int_{Y^*}e^{\lambda\alpha\beta}\sqrt{a}
(s_{\alpha\beta,y}({\bf v})+X_{\alpha\beta})\partial_{\lambda,y}\psi dy-
\frac{1}{2}\int_{Y^*}d^{\alpha\lambda}\sqrt{a}\partial_{\alpha,y}\psi 
\partial_{\lambda,y}\psi dy.$\\
\end{tabular}$$
This yields
$$ \di I({\bf W},\psi)\leq I({\bf W},\Lambda)\leq I({\bf v},\Lambda),$$
for all $({\bf v},\psi)\in {\bf H}^1_{\rm per}(Y^*)/\RR \times  
H^1_{\rm per}(Y^*)/\RR.$\\
Consequenthly, for $\psi=0$ we get
$$I({\bf W},\Lambda)\geq I({\bf W},0)=
\frac{1}{2}\int_{Y^*}
c^{\alpha\beta\lambda\mu}\sqrt{a}\;
(s_{\alpha\beta}({\bf v})+X_{\alpha\beta})
(s_{\lambda\mu,y}({\bf v})+X_{\lambda\mu})dy>0.$$
Moreover, by taking $({\bf v},\psi)=({\bf W},\Lambda)$ in \eq{c4} we obtain
$$\di{\overline c}^{\alpha\beta\tau\theta}X_{\alpha\beta}X_{\tau\theta}=
2I({\bf W},\Lambda)>0.$$
\smallskip

\noindent 
Let us define the function $\Phi:\RR^4\ri\RR$ by
$$\Phi(\xi_{\alpha \beta})=\overline C^{\alpha \beta \gamma \eta} \xi_{\alpha \beta} \xi_{\gamma \eta}.$$
It is to easy that $\Phi$ is continuous in $\RR^4$ endowed with to 
the norm $\parallel \xi \parallel =(\xi_{\alpha \beta}\xi_{\alpha \beta})
^{\frac{1}{2}}.$
Let 
$$B=\{\xi\in\RR^4;\:\xi \mbox{ symmetric },\|\xi\|=1\}.$$
 Since $B$ is compact, $\Phi$ attains its minimum
in $B.$ Then, there exists $c>0$ such that $\Phi\geq c$ in $B,$ that is
$$\di \Phi\left(\frac{\xi_{\alpha\beta}}{\|\xi_{\alpha\beta}\|}\right)\geq c,\quad
\mbox{ for all symmetric tensor }\;\;\xi=(\xi_{\alpha\beta}).$$
This means that $\di {\overline c}^{\alpha\beta\tau\theta}\xi_{\alpha\beta}\xi_{\tau\theta}\geq
c\xi_{\alpha\beta}\xi_{\tau\theta}.$ The coercivity of $( c^{\alpha\beta\tau\theta})$
is now proved.

\noindent (ii) \underline {The ellipticity and symmetry of the tensor}
 $({\overline d}^{\alpha\sigma}).$
\smallskip

First we prove the symmetry. According to \eq{doidouazecisapte} we have
\neweq{d1}
\di {\overline d}^{\alpha\sigma}=
\int_{Y^*}\left[-e_y({\bf z}^\sigma,y_\alpha)+d_y(y_\sigma+\eta^\sigma,y_\alpha)\right]dy,
\endeq
and
\neweq{d2}
\di -\int_{Y^*} e_y(z^\sigma,y_\alpha)dy=
-\int_{Y^*}e_y(z^\sigma,y_\alpha+\eta^\alpha)dy
+\int_{Y^*}e_y(z^\sigma,\eta^\alpha)dy.
\endeq
Now, we take $({\bf v},\psi)=({\bf z}^{\alpha},\eta^{\alpha})$ in 
\eq{doidouazecisapte}. We get
\neweq{d3}
\left\{\begin{tabular}{ll}
$\di-\int_{Y^*}e_y({\bf z}^\sigma,y_\alpha+\eta^\alpha)
dy$&$\di =\int_{Y^*}c_y(z^\sigma,{\bf z}^\alpha)dy,$\\
$\di \int_{Y^*}e_y({\bf z}^\sigma,\eta^\alpha)$&$\di =\int_{Y^*}
d_y(y_\sigma+\eta^{\sigma},\eta^\alpha)dy.$\\
\end{tabular}\right.\endeq
Replacing \eq{d2} and \eq{d3} in \eq{d1} we obtain
$$\di {\overline d}^{\alpha\sigma}=
\int_{Y^*}c_y({\bf z}^\alpha,{\bf z}^\sigma)+\int_{Y^*}d_y(y_\sigma+\eta^\sigma,
y_\alpha+\eta^\alpha)dy. $$
From the above relation, we deduce the symmetry of the tensor 
$({\overline d}^{\alpha\sigma})$.
\medskip

\noindent Let us prove the coercivity of $({\overline d}^{\alpha\sigma})$.
Let $(X_\sigma)$ be a vector.From \eq{doidouazecitrei} we have
\neweq{d4}
\di {\overline d}^{\alpha\sigma}X_\alpha X_\sigma=
-\int_{Y^*}e^{\alpha\lambda\mu}\sqrt{a}\;s_{\lambda\mu,y}({\bf Z}) X_\alpha dy+
\int_{Y^*}d^{\alpha\lambda}\sqrt{a}\;(X_\lambda+\partial_\lambda \Theta)X_\alpha
dy,
\endeq
where ${\bf Z}={\bf z}^\sigma X_\sigma,$ and $\Theta=\eta^\sigma X_\sigma.$
It is easy to see that $({\bf Z},\Theta)$ is the solution of the following 
variational problem
$$\left\{\begin{tabular}{ll}
$\di \int_{Y^*}\left(c_y({\bf Z},{\bf v})+e_y({\bf v},X_\sigma y_\sigma+\Theta)
\right)dy$&$=0,$\\
$\di \int_{Y^*}\left(-e_y({\bf z}^\sigma,\psi)+d_y(y_\sigma+\eta^{\sigma},
\psi)\right)dy$&$=0,$\\
\end{tabular}\right.$$
for all $({\bf v},\psi)\in  {\bf H}^1_{\rm per}(Y^*)/\RR \times 
H^1_{\rm per}(Y^*)/\RR.$ Moreover,
 $({\bf Z},\Theta)$ is a saddle point of the following functional
$$\di J:{\bf H}^1_{\rm per}(Y^*)/\RR \times 
H^1_{\rm per}(Y^*)/\RR\ri\RR,$$
defined by
$$\begin{tabular}{ll}
$\di J({\bf v},\psi)$&
$\di =\frac{1}{2}\int_{Y^*}
c^{\alpha\beta\lambda\mu}\sqrt{a}s_{\alpha\beta,y}({\bf v})s_{\lambda\mu,y}({\bf v})
dy
+\int_{Y^*}e^{\lambda\alpha\beta}\sqrt{a}
s_{\alpha\beta,y}({\bf v})(X_{\lambda}+\partial_{\lambda,y}\psi) dy
$\\
&$\di \;\;\;-
\frac{1}{2}\int_{Y^*}d^{\alpha\lambda}\sqrt{a}(X_\alpha+\partial_{\alpha,y}\psi)
 (X_\lambda+\partial_{\lambda,y}\psi)dy.$\\
\end{tabular}$$
This yields
$$ \di J({\bf Z},\psi)\leq I({\bf Z},\Theta)\leq I({\bf v},\Theta),$$
for all $({\bf v},\psi)\in {\bf H}^1_{\rm per}(Y^*)/\RR \times 
 H^1_{\rm per}(Y^*)/\RR.$\\
By taking $\psi=0$ in the above inequality, we obtain
$$J({\bf Z},\Theta)\geq J({\bf Z},0)=
\frac{1}{2}\int_{Y^*}
d^{\alpha\sigma}\sqrt{a}\;
(X_\alpha+\partial_{\alpha,y}\Theta)
(X_\sigma+\partial_{\sigma,y}\Theta)dy>0.$$
On the other hand, by taking $({\bf v},\psi)=({\bf Z},\Theta)$ in \eq{d4} we obtain
$$\di{\overline d}^{\alpha\sigma}X_{\alpha}X_{\sigma}=
2I({\bf Z},\Theta)>0.$$
With the same proof as for the coercivity of 
$({\overline c}^{\alpha\beta\tau\theta})$
we deduce the existence of $d>0$ such that
$\di {\overline d}^{\alpha\sigma}X_{\alpha}X_{\sigma}\geq
d X_{\alpha}X_{\sigma}.$\\
The uniqueness of $({\bf u},\varphi)$ follows by the
 Lax-Milgram Theorem.

The proof of Theorem \ref {t3} is now complete.$\qed$
\smallskip
\subsection{Corrector result}
We have the follwing convergence 
\begin{equation}
{\cal T}^{\ep}(\gamma_{\alpha\beta}({\bf u}^{\ep}))\we \gamma_{\alpha\beta,x}({\bf u})+s_{\alpha\beta,y}({\bf u}^1)\quad\mbox{ weakly in }\,{\bf L}^2(\omega\times Y^*),
\label{con1}
\end{equation}
\begin{equation}
{\cal T}^{\ep}(\nabla\varphi^{\ep})\we \nabla_{x}\varphi+\nabla_{y}\varphi^1\quad\mbox{ weakly in }\,L^2(\omega\times Y^*).
\label{con2}
\end{equation}
The convergence of energies follow easily fro the above relations. Moreover, 
the weak convergences in (\ref {con1}) and (\ref {con2}) are actually strong
\begin{equation}
{\cal T}^{\ep}(\gamma_{\alpha\beta,x}({\bf u}^{\ep}))-\gamma_{\alpha\beta,x}({\bf u})-s_{\alpha\beta,y}({\bf u}^1)\ri0\quad\mbox{ strongly in }\,{\bf L}^2(\omega\times Y^*),
\label{sto1}
\end{equation}
\begin{equation}
{\cal T}^{\ep}(\nabla_x \varphi^{\ep})-\nabla_{x}\varphi-\nabla_{y}\varphi^1\ri 0\quad\mbox{ strongly in }\,L^2(\omega\times Y^*).
\label{sto2}
\end{equation}
Now, we can state the following corrector result:
\begin{teo}$(${\bf correctors}$)$. One has the following strong convergences :
\begin{equation}
\gamma_{\alpha\beta,x}({\bf u}^{\ep})-\gamma_{\alpha\beta,x}({\bf u})-{\cal U}^{\ep}(s_{\alpha\beta,y}({\bf u}^1))\ri0\quad\mbox{ strongly in }\,
{\bf L}^2(\omega),
\label{tof1}
\end{equation}
\begin{equation}
\nabla_x \varphi^{\ep}-\nabla_{x}\varphi-{\cal U}^{\ep}(\nabla_{y}\varphi^1)\ri 0\quad\mbox{ strongly in }\,L^2(\omega).
\label{tof2}
\end{equation}
\end{teo}
\noindent{\bf Proof.} Using convergences (\ref {sto1})-(\ref {sto2}) and Theorem \ref{topu} we have
\neweq{cdem1}
\gamma_{\alpha\beta,x}({\bf u}^{\ep})-{\cal U}^{\ep}(\gamma_{\alpha\beta,x}({\bf u}))-{\cal U}^{\ep}(s_{\alpha\beta,y}({\bf u}^1))\ri0\quad\mbox{ strongly in }\,{\bf L}^2(\omega),
\endeq
\neweq{cdem2}
\nabla_x \varphi^{\ep}-{\cal U}^{\ep}(\nabla_{x}\varphi)-{\cal U}^{\ep}(s_{\alpha\beta,y}({\bf u}^1))\ri0\quad\mbox{ strongly in }\,L^2(\omega).
\endeq
But $\gamma_{\alpha\beta,x}({\bf u})\in {\bf L}^2(\omega)$ and $\nabla_x \varphi\in L^2(\omega)$ 
so, by Proposition \ref{p3} (i) we get 
\neweq{cdem3}
{\cal U}^{\ep}(\gamma_{\alpha\beta,x}({\bf u}))\ri\gamma_{\alpha\beta,x}({\bf u})\quad\mbox{ strongly in }\,{\bf L}^2(\omega),
\endeq
\neweq{cdem4}
{\cal U}^{\ep}(\nabla_x \varphi)\ri\nabla_x \varphi\quad\mbox{ strongly in }\,L^2(\omega).
\endeq
From \eq{cdem1}-\eq{cdem4} we deduce the convergence \eq{tof1} and \eq{tof2}. $\qed$

\section{The bending problem}
\subsection{The convergence results}
We now consider the variational bending problem
\neweq{cunu}
\begin{tabular}{ll}
$\di \frac{2}{3}\int_{\omega^{\ep}}C^{\alpha\beta\lambda\mu,\ep}\Upsilon_{\alpha\beta}({\bf u}^{\ep})\Upsilon_{\lambda\mu}({\bf v})\sqrt{a^{\ep}}dx$\\
$\di\qquad=\int_{\omega^{\ep}}\left(\int_{-1}^1f^{i}(x_1,x_2,z)dz\right)v^i\sqrt{a^{\ep}}~dx+\int_{\Gamma_{+}^{\ep}\cup\Gamma_{-}^{\ep}}q^{i}v^i\sqrt{a^{\ep}}d\Gamma,$
\end{tabular}\endeq
for all ${\bf v} \in {\bf W}(\omega^{\ep})$, where
$$C^{\alpha\beta\lambda\mu,\ep}=C^{\alpha\beta\lambda\mu}(\frac{x}{\ep}),~~a^{\ep}=a(\frac{x}{\ep})$$
$$$$
$${\bf W}(\omega^{\ep})=\Big\{{\bf w}\in H^1(\omega^{\ep})\times H^1(\omega^{\ep})\times H^2(\omega^{\ep});\;\gamma_{\alpha\beta}({\bf w})=0\mbox{ and } w^i=\partial_{\nu}w^3=0~\mbox{ on }\gamma^0_{\ep}\subset \partial\omega^{\ep}\Big\}  $$
We have
$${\bf f}=(f^i)\in {\bf L}^2(\omega^{\ep}\times [-1,1]),\quad 
{\bf q}=(q^i)\in{\bf L}^2(\Gamma_{+}^{\ep}\cup\Gamma_{-}^{\ep}),$$

\neweq{ctrei}
\di \Upsilon_{\alpha\beta}({\bf v})=-\Big(\partial_{\alpha\beta}^2v_3-v_{\rho}\Big(-\partial_\beta b^{\rho}_\alpha+b_\beta^\gamma \Gamma^\rho_{\alpha\gamma}+\gamma^\delta_{\alpha\beta}b^\rho_\delta\Big)-c_{\alpha\beta}v_3+b^\rho_\alpha\partial_\beta v_\nu-\Gamma_{\alpha\beta}^\delta \partial_\delta v_3\Big).
\endeq
Due to the ellipticity and symmetry of the bending tensor, by using the Lax-Milgram Theorem, we can deduce the existence and uniqueness to solution ${\bf u}^{\ep}$ of the variational problem \eq{cunu}.
\smallskip
\\
Let us denote 
$$\di F^i(x_1,x_2)=\int_{-1}^1f^{i}(x_1,x_2,z)~dz.$$ 
Thus, the equation \eq{cunu} takes the form
\neweq{csase}
\di \frac{2}{3}\int_{\omega^{\ep}}C^{\alpha\beta\lambda\mu,\ep}\sqrt{a^{\ep}}\Upsilon_{\alpha\beta}({\bf u}^{\ep})\Upsilon_{\rho\sigma}({\bf v})dx=
\int_{\omega^{\ep}}F^iv^i\sqrt{a^{\ep}}~dx+\int_{\Gamma_{+}^{\ep}\cup\Gamma_{-}^{\ep}}h^{i}v^i\sqrt{a^{\ep}}~d\Gamma,
\endeq
for all $\;{\bf v}\in {\bf W}(\omega^{\ep}).$
\begin{teo}\label{t4}
The sequence $\left({\cal T}^{\ep}({\bf u}^{\ep})\right)_{\ep}$ weakly converges to ${\bf u} \in {\bf W}(\omega)$ which is the unique solution of the homogenized problem 
\neweq{csapte}
\di\frac{2}{3}\int_{\omega}\overline C^{\alpha\beta\rho\sigma}\,\Upsilon_{\alpha\beta}({\bf u})\,\Upsilon_{\rho\sigma}({\bf v})~dx=\mid Y^* \mid_a\int_{\omega}F^i v^i~dx+\mid Y^* \mid_a\int_{\Gamma_{+}\cup\Gamma_{-}}h^{i}v^i~d\Gamma,
\endeq
for all ${\bf v}\in {\bf W}(\omega),$ where
$$\mid Y^* \mid_a=\int_{Y^*}\sqrt{a(y)}~dy,$$
\neweq{copt}
\di \overline C^{\alpha\beta\rho\sigma}=\int_{Y^*}C^{\tau\theta\rho\sigma}\sqrt{a}\,\left[\delta_{\alpha\tau}\delta_{\beta\theta}+\partial_{\tau\theta,y}^2 {\bf w}^{\alpha\beta}\right]dy,
\endeq
and the local functions $w^{\tau\theta}$ verifes the local problems 
\neweq{cnoua}
\left\{\begin{tabular}{ll}
$\di\frac{\partial}{\partial y_\rho}\left\{C^{\alpha\beta\rho\sigma}\sqrt{a}\,\left[\delta_{\alpha\tau}\delta_{\beta\theta}+\partial_{\alpha\beta,y}^2{\bf w}^{\tau\theta}\right]\right\}=0$&$\quad\mbox{ in }\;\omega\times Y^*,$\\
${\bf w}^{\tau\theta}\quad Y^*\mbox{-periodic}.$\end{tabular}\right.
\endeq
\end{teo}
\noindent{\bf Remark.} We can give another expression for the homogenized bending coefficients, that as form :
$$\di \overline C^{\alpha\beta\rho\sigma}=
\int_{Y^*} C^{\tau\theta\rho\sigma}\sqrt{a}
\partial^2_{\tau\theta,y}[\Pi^{\alpha\beta}+{\bf w}^{\alpha\beta}]~dy,$$
where $\Pi^{\alpha\beta}=\frac{1}{2}y_\alpha y_\beta.$\\
\noindent{\bf Proof.} To obtain an {\it a priori estimate} for ${\bf u}^{\ep},$ we choose 
${\bf v}={\bf u}^{\ep}$ in \eq{csase}. It follows
$$\di \int_{\omega}C^{\alpha\beta\rho\sigma,\ep}\sqrt{a^\ep}\Upsilon_{\alpha\beta}({\bf u}^{\ep})\Upsilon_{\rho\sigma}({\bf u}^{\ep})~dx\leq M\|{\bf u}^{\ep}\|_{L^2(\omega_{\ep})},$$
where $M>0$ is a positive constant that not depend on ${\ep}.$
Using the ellipticity of the tensor $C^{\alpha\beta\rho\sigma,\ep}$ we get
$$ \di c\int_{\omega^{\ep}}\Upsilon_{\alpha\beta}^2({\bf u}^{\ep})dx\leq M\|{\bf u}^{\ep}\|_{L^2(\omega^{\ep})}.$$
Now, the Korn and Poincar\'e's inequalities in perforated domains imply
$$\|{\bf u}^{\ep}\|_{{\bf W}(\omega^{\ep})}\leq C,$$
where $C$ does not depend on $\ep.$ Then, up to a subsequence, $({\bf u}^{\ep})$ weakly converges to some ${\bf u}\in {\bf W}(\omega).$ By Theorem 1.2 we deduce that there exists ${\bf u}^2\in L^2(\omega;H^2_{\rm per}(Y^*)/{\RR})$ such that
\neweq{czece}
{\cal T}^{\ep}({\bf u}^{\ep})\we {\bf u}\quad\mbox{ weakly in }\,{\bf L}^2(\omega\times Y^*),
\endeq
\neweq{cunspe}
{\cal T}^{\ep}(\nabla {\bf u}^{\ep})\we \nabla {\bf u}\quad\mbox{ weakly in }\,{\bf L}^2(\omega\times Y^*),
\endeq
and
\neweq{cdoispe}
{\cal T}^{\ep}(\nabla^2 {\bf u}^{\ep})\we \nabla^2{\bf u}+\nabla^2_y {\bf u}^2\quad\mbox{ weakly in }\,L^2(\omega\times Y^*).
\endeq
The linearity of ${\cal T}^{\ep}$ implies
\neweq{ctreispe}
{\cal T}^{\ep}(\Upsilon_{\alpha\beta}({\bf u}^{\ep}))\we \Upsilon_{\alpha\beta}({\bf u})+\partial_{\alpha\beta,y}^2{\bf u}^2\quad\mbox{ weakly in }\,L^2(\omega\times Y^*).
\endeq
Using now the properties of the unfolding operator ${\cal T}^{\ep}$, in \eq{csase}, we get
$$\begin{tabular}{ll}
$\di \frac{2}{3}\int_{\omega^{\ep}\times Y^*}{\cal T}^{\ep}
(C^{\alpha\beta\rho\sigma,\ep}\sqrt{a^{\ep}})\,
{\cal T}^{\ep}(\Upsilon_{\alpha\beta}
({\bf u}^{\ep}))\,{\cal T}^{\ep}(\Upsilon_{\rho\sigma}({\bf v}))
dxdy$\\
$\qquad\qquad\qquad
\hspace{2cm}\di=\int_{\omega^{\ep}\times Y^*}{\cal T}^{\ep}(F^iv^i\sqrt{a^{\ep}})~dxdy+\int_{(\Gamma_{+}^{\ep}\cup\Gamma^{\ep}_{-})\times Y^*}{\cal T}^{\ep}(h^{i}v^i\sqrt{a^{\ep}})~d\Gamma dy,$\\
\end{tabular}$$
that is
\neweq{cpaispe}
\begin{tabular}{ll}
$\di\frac{2}{3}\int_{\omega^{\ep}\times Y^*}C^{\alpha\beta\rho\sigma}
\sqrt{a}\,{\cal T}^{\ep}(\Upsilon_{\alpha\beta}({\bf u}^{\ep}))\,
{\cal T}^{\ep}(\Upsilon_{\rho\sigma}({\bf v}))dxdy$\\
$\qquad\qquad
\hspace{2cm}\di=\int_{\omega^{\ep}\times Y^*}{\cal T}^{\ep}
(F^iv^i\sqrt{a^{\ep}})~dxdy+
\int_{(\Gamma_{+}^{\ep}\cup\Gamma^{\ep}_{-})\times Y^*}
{\cal T}^{\ep}(h^{i}v^i\sqrt{a^{\ep}})~d\Gamma dy,$\\
\end{tabular}
\endeq
We chose as a test function in \eq{cpaispe}
$$\di {\bf v}^\ep(x)={\bf v}_1(x)+\ep^2{\bf v}_2\left(x,\frac{x}{\ep}\right),
\quad{\bf v}_1\in{\cal D}(\omega);~{\bf v}_2\in{\cal D}(\omega;C^{\infty}_{\rm per}(Y^*)).$$
Then
\begin{eqnarray}
\di \nabla_x{\bf v}^\ep(x)&=&\nabla{\bf v}_1(x)+\ep^2\nabla_x{\bf v}_2\left(x,\frac{x}{\ep}\right)+\ep\nabla_y{\bf v}_2\left(x,\frac{x}{\ep}\right),\nonumber \\
\di \nabla_x^2{\bf v}^\ep(x)&=&\nabla^2{\bf v}_1(x)+\ep^2\nabla^2_x{\bf v}_2\left(x,\frac{x}{\ep}\right)+2\ep\nabla_x\nabla_y{\bf v}_2\left(x,\frac{x}{\ep}\right)+\nabla^2_y{\bf v}^2\left(x,\frac{x}{\ep}\right).
\end{eqnarray}
It follows that
\begin{equation}
\begin{array}{cllll} \left\{ \begin{array}{cllll}
{\bf v}^\ep(x)&\ri&{\bf v}_1(x)\quad&\mbox{ strongly in }\,{\bf L}^2(\omega),\\
\nabla_x{\bf v}^\ep(x)&\ri&\nabla{\bf v}_1(x)\quad&\mbox{ weakly in }\,{\bf L}^2(\omega),\\
\nabla^2_x{\bf v}^\ep(x)&\we&\nabla^2{\bf v}_1(x)+\nabla^2_y{\bf v}_2\left(x,y\right)\quad&\mbox{ weakly in }\,{\bf L}^2(\omega\times Y^*).
\end{array} \right.
\end{array}
\label{ccincispe}
\end{equation}
And
\neweq{csaispe}
{\cal T}^{\ep}(\Upsilon_{\rho\sigma}({\bf v}^\ep))\we \Upsilon_{\rho\sigma}({\bf v}_1)+\partial_{\rho\sigma,y}^2{\bf v}_2\quad\mbox{ weakly in }\,{\bf L}^2(\omega\times Y^*).
\endeq
Using to \eq{ccincispe} and \eq{csaispe}, and passing to the limit in \eq{cpaispe} with $\ep \searrow 0,$ we get
\neweq{csaptispe}
\begin{tabular}{ll}
$ \di \frac{2}{3}\int_{\omega\times Y^*}C^{\alpha\beta\rho\sigma}\sqrt{a}\,\Big[\Upsilon_{\alpha\beta}({\bf u})+\partial_{\alpha\beta,y}^2{\bf u}^2\Big]\,\Big[\Upsilon_{\rho\sigma}({\bf v}_1)+\partial_{\rho\sigma,y}^2 {\bf v}_2\Big]dxdy$\\
$\qquad\qquad\qquad\hspace{3cm}
\di=\mid Y^* \mid_a\int_{\omega}F^i v_1^i~dx+\mid Y^* \mid_a\int_{\Gamma_{+}\cup\Gamma_{-}}h^{i}v_1^i~d\Gamma.$\\
\end{tabular}
\endeq
We now let ${\bf v}_2=0$ in \eq{csaptispe}. We obtain
\neweq{coptispe}
\begin{tabular}{ll}
$ \di\frac{2}{3}\int_{\omega\times Y^*}C^{\alpha\beta\rho\sigma}\sqrt{a}\,\Big[\Upsilon_{\alpha\beta}({\bf u})+\partial_{\alpha\beta,y}^2{\bf u}^2\Big]\,\Upsilon_{\rho\sigma}({\bf v}_1)~dxdy$\\$\qquad\qquad\qquad\hspace{3cm}
\di=\mid Y^* \mid_a\int_{\omega}F^i v_1^i~dx+\mid Y^* \mid_a\int_{\Gamma_{+}\cup\Gamma_{-}}h^{i}v_1^i~d\Gamma.$\\
\end{tabular}
\endeq
By density, we deduce
\neweq{cnouaspe}
\begin{tabular}{ll}
$ \di\frac{2}{3}\int_{\omega\times Y^*}C^{\alpha\beta\rho\sigma}\sqrt{a}\,\Big[\Upsilon_{\alpha\beta}({\bf u})+\partial_{\alpha\beta,y}{\bf u}^2\Big]\,\Upsilon_{\rho\sigma}({\bf v})~dxdy$\\
$\qquad\qquad\qquad\hspace{3cm}
\di=\mid Y^* \mid_a\int_{\omega}F^iv_1^i~dx+\mid Y^* \mid_a\int_{\Gamma_{+}\cup\Gamma_{-}}h^{i}v_1^i~d\Gamma,$\\
\end{tabular}
\endeq
for all ${\bf v}\in L^2(\omega;H^2_{\rm per}(Y^*)/\RR).$ In what follows we chose 
\neweq{cdouazeci}
\di {\bf u}^2=\Upsilon_{\tau\theta}({\bf u}){\bf w}^{\tau\theta},\quad {\bf w}^{\tau\theta}\in L^2(\omega;H^2_{\rm per}(Y^*)/\RR).\endeq
Then
\neweq{cdoiunu}
\di \partial_{\alpha\beta,y}{\bf u}^2=\Upsilon_{\tau\theta}({\bf u})\partial_{\alpha\beta,y}^2{\bf w}^{\tau\theta}.
\endeq
Replacing \eq{cdoiunu} in \eq{cnouaspe} we obtain
\neweq{cdoidoi}
\begin{tabular}{ll}
$ \di \frac{2}{3}\int_{\omega\times Y^*}C^{\alpha\beta\rho\sigma}\sqrt{a}\,\Big[\Upsilon_{\alpha\beta}({\bf u})+\Upsilon_{\tau\theta}({\bf u})\partial_{\alpha\beta,y}^2{\bf w}^{\tau\theta}\Big]\,\Upsilon_{\rho\sigma}({\bf v})~dxdy$\\
$\qquad\qquad\qquad\hspace{3cm}
\di=\mid Y^* \mid_a\int_{\omega}F^i v^i~dx+\mid Y^* \mid_a\int_{\Gamma_{+}\cup\Gamma_{-}}h^{i} v^i~d\Gamma,$\\
\end{tabular}
\endeq
for all ${\bf v}\in L^2(\omega;H^2_{\rm per}(Y^*)/\RR).$ It follows that
$$\begin{tabular}{ll}
$ \di\frac{2}{3}\int_{\omega\times Y^*}C^{\tau\theta\rho\sigma}\sqrt{a}\,\Big[\Upsilon_{\tau\theta}({\bf u})+\Upsilon_{\alpha\beta}({\bf u})\partial_{\tau\theta,y}^2{\bf w}^{\alpha\beta}\Big]\,\Upsilon_{\rho\sigma}({\bf v})~dxdy$\\
$\qquad\qquad\qquad\hspace{3cm}
\di=\mid Y^* \mid_a\int_{\omega}F^iv^i~dx+\mid Y^* \mid_a\int_{\Gamma_{+}\cup\Gamma_{-}}h^{i}v^i~d\Gamma,$\\
\end{tabular}$$
for all ${\bf v}\in L^2(\omega;H^2_{\rm per}(Y^*)/\RR) .$ This yields
\neweq{cdoitrei}
\begin{tabular}{ll}
$ \di\frac{2}{3}\int_{\omega}\di\Big\{\int_{Y^*} C^{\tau\theta\rho\sigma}\sqrt{a}\Big[\delta_{\alpha\tau}\delta_{\beta\theta}+\partial_{\tau\theta,y}^2{\bf w}^{\alpha\beta}\Big]dy\Big\}\Upsilon_{\alpha\beta}({\bf u})\,\Upsilon_{\rho\sigma}({\bf v})~dx$\\
$\qquad\qquad\qquad\hspace{3cm}
\di=\mid Y^* \mid_a\int_{\omega}F^iv^i~dx+\mid Y^* \mid_a\int_{\Gamma_{+}\cup\Gamma_{-}}h^{i}v^i~d\Gamma,$\\
\end{tabular}\endeq
for all ${\bf v}\in L^2(\omega;H^2_{\rm per}(Y^*)/\RR).$ \\
If we denote
\neweq{cdoipatru}
\di \overline C^{\alpha\beta\rho\sigma}=\int_{Y^*} C^{\tau\theta\rho\sigma}\sqrt{a(y)}\,\Big[\delta_{\alpha\tau}\delta_{\beta\theta}+\partial_{\tau\theta,y}^2{\bf w}^{\alpha\beta}\Big]dy,
\endeq
by \eq{cdoitrei} we get the homogenized equation which corresponds to \eq{csase}:
\neweq{cdoicinci}
\di\frac{2}{3}\int_{\omega}\overline C^{\alpha\beta\rho\sigma}\,\Upsilon_{\alpha\beta}({\bf u})\,\Upsilon_{\rho\sigma}({\bf v})dx=\mid Y^* \mid_a\int_{\omega}F^iv^i~dx+\mid Y^* \mid_a\int_{\Gamma_{+}\cup\Gamma_{-}}h^{i}v^i~d\Gamma,
\endeq
for all ${\bf v} \in L^2(\omega;H^2_{\rm per}(Y^*)/\RR).$

Let us find now the equations verified by the local functions ${\bf w}^{\tau\theta}.$
\\
Letting $\phi=0$ in \eq{csaptispe} we get
$$\di\int_{\omega\times Y^*}C^{\alpha\beta\rho\sigma}\sqrt{a}\,\Big
[\Upsilon_{\alpha\beta}({\bf u})+\partial_{\alpha\beta,y}^2{\bf u}^2\Big]
\partial_{\rho\sigma,y}^2{\bf v}_2(x,y)~dxdy=0.$$
By density it follows that
\neweq{cdoisase}
\di\int_{\omega\times Y^*}C^{\alpha\beta\rho\sigma}(y)\sqrt{a}\,\Big[\Upsilon_{\alpha\beta}({\bf u})+\partial_{\alpha\beta,y}{\bf u}^2\Big]\,
\partial_{\rho\sigma,y}^2{\bf v}(x,y)~dxdy=0,\endeq
for all ${\bf v}\in {\bf L}^2(\omega;H^2_{\rm per}(Y^*)/\RR).$\\
Using now \eq{cdoiunu} in \eq{cdoisase} we have
$$\di\frac{\partial}{\partial y_\rho}
\Big\{\Upsilon_{\alpha\beta}({\bf u})C^{\alpha\beta\rho\sigma}\sqrt{a}\,\Big[\delta_{\alpha\tau}\delta_{\beta\theta}+\partial_{\alpha\beta,y}^2{\bf w}^{\tau\theta}\Big]\Big\}=0\quad\mbox{ in }\;\omega\times Y^*$$
and so
\neweq{cdoisapte}
\di\frac{\partial}{\partial y_\rho}\left\{C^{\alpha\beta\rho\sigma}\sqrt{a}\,\left[\delta_{\alpha\tau}\delta_{\beta\theta}+\partial_{\alpha\beta,y}^2{\bf w}^{\tau\theta}\right]\right\}=0\quad\mbox{ in }\;\omega\times Y^*.
\endeq
In order to establish the existence and uniqueness of the solution of (\ref {csapte}), if suffices to prove the coercivity of $\overline C^{\alpha\beta\gamma\theta}$ in the following sense
$$\exists \Lambda_C \neq \Lambda_C (\varepsilon) >0,~~\forall (\xi_{\alpha\beta})_{\alpha\beta}~:~\xi_{\alpha\beta}=\xi_{\beta\alpha},~~\overline C_{\alpha \beta\gamma\eta} \xi_{\alpha\beta} \xi_{\gamma\eta}\geq  \Lambda_C \xi_{\alpha\beta } \xi_{\alpha\beta}$$

\noindent{\it Symmetry}\\
It is easy to check that 
$$\overline C^{\alpha \beta \gamma \eta}=\overline C^{\beta \alpha \gamma \eta}=\overline C^{\alpha \beta \eta \gamma}$$
It suffices to prove  $$\overline C^{\alpha \beta \gamma \eta}=\overline C^{\gamma \eta\alpha\beta}$$ 
Starting from the definition (\ref {cdoipatru}) of the $\mathcal{\overline C}=(\overline C^{\alpha \beta \gamma \eta})$, the homogenized bendig tensor is evaluated by
\begin{eqnarray}
\displaystyle \overline C^{\alpha \beta \gamma \eta}&=
&\int_{Y^*} C^{\alpha\beta\delta\tau }\sqrt{a}\;
\partial^2_{\partial \delta\tau,y}[\Pi^{\gamma\eta}+{\bf w}^{\gamma\eta}]~dy\nonumber \\
&=&\int_{Y^*} C^{\zeta\varsigma\delta\tau }\sqrt{a}\;
\partial^2_{\partial \delta\tau,y}[\Pi^{\gamma\eta}+{\bf w}^{\gamma\eta}]
\delta_{\alpha\zeta}\delta_{\beta\varsigma}~dy\nonumber \\
&=&\int_{Y^*} C^{\zeta\varsigma\delta\tau }\sqrt{a}\;
\partial^2_{\partial \delta\tau,y}[\Pi^{\gamma\eta}+{\bf w}^{\gamma\eta}]\;
\partial^2_{\partial \zeta\varsigma,y}\Pi^{\alpha\beta}~dy\nonumber \\
&=&\int_{Y^*} C^{\zeta\varsigma\delta\tau }\sqrt{a}\;
\partial^2_{\partial \delta\tau,y}[\Pi^{\gamma\eta}+{\bf w}^{\gamma\eta}]\;
\partial^2_{\partial \zeta\varsigma,y}
[\Pi^{\alpha\beta}+{\bf w}^{\alpha\beta}]~dy\nonumber \\
&-&\int_{Y^*} C^{\zeta\varsigma\delta\tau }\sqrt{a}\;
\partial^2_{\partial \delta\tau,y}[\Pi^{\gamma\eta}+{\bf w}^{\gamma\eta}]
\partial^2_{\partial \zeta\varsigma,y}
{\bf w}^{\alpha\beta}~dy.
\end{eqnarray}
By multiplying (\ref {cnoua}) by ${\bf w}^{\alpha\beta}$ and integrating by parts, we prove :
$$
\int_{Y^*} C^{\zeta\varsigma\delta\tau }\sqrt{a}\;
\partial^2_{\partial \delta\tau,y}[\Pi^{\gamma\eta}+{\bf w}^{\gamma\eta}]
\partial^2_{\partial \zeta\varsigma,y}
{\bf w}^{\alpha\beta}~dy=0.
$$
It follows that
\begin{equation}
\overline C^{\alpha \beta \gamma \eta}=\int_{Y^*}C^{\delta \eta \zeta \nu}
\sqrt{a}\;
\partial^2_{\partial \delta\tau,y}\Big (\Pi^{\gamma\eta}+{\bf w}^{\gamma\eta}
ig )
\partial^2_{\partial \zeta\varsigma,y}
\Big (\Pi^{\gamma \eta}+{\bf w}^{\gamma \eta}(y)\Big )~dy.
\label{toj}
\end{equation}
From (\ref {toj}), we deduce $\overline C^{\alpha\beta\gamma\theta}=\overline C^{\gamma\theta\alpha\beta}$.

\noindent{\it Ellipticity}\\
Let $(\xi_{\alpha\beta})_{\alpha\beta}$ be a symmetric tensor
($\xi_{\alpha\beta}=\xi_{\beta\alpha}$) and set
$$ \tau_{\delta\eta}=\xi_{\alpha\beta}\partial^2_{\delta\eta,y}
\Big (\Pi^{\alpha \beta}+{\bf w}^{\alpha \beta}\Big ).$$
Using now the coercivity of tensor $C^{\lambda \mu \iota \varsigma}(x,y)$ 
and the fact that $a\not\equiv 0,$ we  can write
\begin{equation}
\overline C^{\alpha \beta \gamma \eta} \xi_{\alpha \beta} \xi_{\gamma \eta} \geq \int_{Y^*} C^{\delta \eta \zeta \nu}\sqrt{a}\tau_{\delta\eta} \tau_{\zeta\nu} ~dy \geq c \int_{Y^*}\tau_{\delta\eta} \tau_{\delta\eta}~dy.
\label{eeli}
\end{equation}
We claim that the second integral in (\ref {eeli}) is positive. Assume 
the contrary. Then
\begin{equation}
\forall~~(\delta,\eta) \in \{1,2\}^2,~~~\tau_{\delta\eta}=
\xi_{\alpha\beta}\partial^2_{\delta\eta,y}\Big (\Pi^{\alpha \beta}-{\bf w}^{\alpha \beta}\Big )=0. 
\label{absurde}
\end{equation}
It follows that
$$ \partial^2_{\delta\eta,y}\Big (
\xi_{\alpha\beta}(\Pi^{\alpha \beta}-{\bf w}^{\alpha \beta})\Big )=0. $$
This implies that
$$\xi_{\alpha\beta}(\Pi^{\alpha \beta}-{\bf w}^{\alpha \beta})\Big )=
a_{\iota} y_{\iota}+b,$$
for some constants $a_{\iota}$ and $b,$ $\iota=1,2.$
This yields 
$${\bf w}^{\alpha\beta}\xi_{\alpha\beta}=\Pi^{\alpha\beta}\xi_{\alpha\beta}+a_{\iota}y_{\iota}+b.$$
Since $\xi\not\equiv 0,$ we can find an indice $(\alpha,\beta)$ 
such that $\xi_{\alpha\beta}\neq 0.$ In this case, the left-hand 
side of the above equality is $Y^*$-periodic, but the right-hand side is not,
this is clearly a  contradiction.
Then the second integral of (\ref {eeli}) is positive and so
$$\overline C^{\alpha \beta \gamma \eta} \xi_{\alpha \beta} \xi_{\gamma \eta} >0 ~~~~\forall~(\xi_{\alpha\beta})\neq 0~\mbox{symmetric }.$$
Let us define the function $\Psi:\RR^4\ri\RR$ by
$$\Psi(\xi_{\alpha \beta})=\overline C^{\alpha \beta \gamma \eta} \xi_{\alpha \beta} \xi_{\gamma \eta}.$$
It is easy to see that $\Psi$ is continuous in $\RR^4$ endowed with to 
the norm 
$$\parallel \tau \parallel =(\tau_{\alpha \beta}\tau_{\alpha \beta})
^{\frac{1}{2}}.$$
Since $\Phi$ attains its minimum on the unit sphere in $\RR^4$ and 
$\Psi>0$ for all symmetric tensor $(\xi_{\alpha\beta})\not\equiv 0,$
we can conclude that there exists
$M>0$ such that 
$$\Psi(\frac{\xi_{\alpha \beta}}{\parallel \xi \parallel}) \geq M,
~~~\mbox{ for all symmetric tensor }~ 
(\xi_{\alpha\beta})\not\equiv 0.$$
From the above inequality we deduce
$$\overline C^{\alpha \beta\gamma\eta} 
\xi_{\alpha\beta} \xi_{\gamma\eta}\geq 
M \xi_{\alpha\beta } \xi_{\alpha\beta}.$$
The uniqueness of the solution of \eq{csapte} follows now by using the 
Lax-Milgram Theorem.$\qed$
\subsection{Corrector result}
We have the follwing convergence :
\begin{equation}
{\cal T}^{\ep}(\Upsilon_{\rho\sigma}({\bf u}^\ep))\we \Upsilon_{\rho\sigma,x}({\bf u})+\partial_{\rho\sigma,y}^2{\bf u}_2\quad\mbox{ weakly in }\,{\bf L}^2(\omega\times Y^*).
\label{con}
\end{equation}
The convergence of energies is also proved easily, and implies in particular that the weak convergences in (\ref {con}) is actually strong
\begin{equation}
{\cal T}^{\ep}(\Upsilon_{\rho\sigma,x}({\bf u}^\ep))\ri \Upsilon_{\rho\sigma,x}({\bf u})+\partial_{\rho\sigma,y}^2{\bf u}_2\quad\mbox{ strongly in }\,{\bf L}^2(\omega\times Y^*).
\label{conn}
\end{equation}
\begin{teo}$(${\bf correctors}$)$. One has the following strong convergence :
$$\Upsilon_{\rho\sigma,x}({\bf u}^\ep)-\Upsilon_{\rho\sigma,x}({\bf u})-{\cal U}^{\ep}(\partial_{\rho\sigma,y}^2{\bf u}_2)\ri 0
\quad\mbox{ strongly in }\,{\bf L}^2(\omega).$$
\end{teo}
\noindent{\bf Proof.} We have already seen (see (\ref {conn})) that
$$
{\cal T}^{\ep}(\Upsilon_{\rho\sigma,x}({\bf u}^\ep))-\Upsilon_{\rho\sigma,x}({\bf u})-\partial_{\rho\sigma,y}^2{\bf u}_2 \ri 0\quad\mbox{ strongly in }\,{\bf L}^2(\omega\times Y^*),
$$
which, by Theorem \ref{topu} is equivalent to
$$
\Upsilon_{\rho\sigma,x}({\bf u}^\ep)-{\cal U}^{\ep}(\Upsilon_{\rho\sigma,x}({\bf u})+\partial_{\rho\sigma,y}^2{\bf u}_2)\ri 0\quad\mbox{ strongly in }\,{\bf L}^2(\omega).
$$
But $\Upsilon_{\rho\sigma,x}({\bf u})\in {\bf L}^2(\omega)$, so from (i) of Proposition \ref{p3} on has ${\cal U}^{\ep}(\Upsilon_{\rho\sigma,x}({\bf u}))\ri \Upsilon_{\rho\sigma,x}({\bf u})$ strongly in ${\bf L}^2(\omega)$, whence the desired result.$\qed$
\section{Conclusion}
In this paper we have rigously established the limiting equations modelling the behavior of a thin piezoelectric perforated shells, i.e., we have explicity described the {\it homogenized coefficients} of the elastic, dielectric and coupling tensors (for details, see \cite{Mech2}).\\\\
{\bf Acknowledgment.} This work has been supported in part by the Ministry for higher education and scientific research of Algeria (University of Oran, Departement of Mathematics). The author is grateful to Professor Bernadette Miara for helpful discussions.


\begin{thebibliography}{16} \frenchspacing

\bibitem{Ben} Bensoussan  A., Lions J.L., Papanicolaou G. {\it Asymptotic Analysis for Periodic Structures}, North  Holland, Amsterdam (1978).
\bibitem{Cior} Cioranescu D., Damlamian A., Griso G. {\it Periodic unfolding and homogenization}, C. R. Acad. Sci. Paris, Ser. I {\bf 334} (2002) 99-104.
\bibitem{Gher1} Ghergu M., Griso G., Mechkour H., Miara B. {\it Homog\'en\'eisation de coques minces pi\'ezo\'electriques perfor\'ees}. C. R. Acad. Sci. Paris, Ser. II : M\'ecanique {\bf 333} (2005) 249-255.
\bibitem{Gher3} Ghergu M., Griso G., Labat B., Mechkour H., Miara B., Rohan E., Zidi M. {\it Homog\'en\'eisation et pi\'ezo\'electricit\'e. Aide \`a la conception d'un bio-mat\'eriau}. Annals of University of Craoiva. Math. Comp. Sci. Ser. {\bf 32} (2005): 9-15. 
\bibitem{Haen} Haenel Ch. {\it Analyse et simulation num\'erique de coques pi\'ezo\'electriques}. PhD thesis, University of Paris 6 (2000).
\bibitem{LewTel} Lewinski T., Telega J.J. {\it Plates, Laminates and shells: Asymptotic analysis and homogenization.} Ser. Adv. Math. Applied Sciences - Vol. 52, World Scientific (2000).
\bibitem{Mech0} Mechkour H. {\it Two-scale homogenization of periodic perforated piezoelectric structures.} Preprint, 09/2007 (http://arxiv.org/abs/0709.1079).
\bibitem{Mech01} Mechkour H. \emph{Comportement macroscopique d'une plaque pi\'ezo\'electrique p\'eriodiquement perfor\'ee}, In Proc. of The 3rd Conference on Trends in Applied Mathemtics in Tunisia, Algeria, Morocco (TAMTAM-3), 16-18 April 2007, Alger, Algeria (in french).
\bibitem{Mech2} Mechkour H {\it Homog\'en\'eisation et simulation num\'erique de structures pi\'ezo\'electriques perfor\'ees et lamin\'ees}. PhD thesis, University of Marne-La-Vall\'ee 2004 (in french).
\bibitem{Mech1} Mechkour H., Miara B. \emph{ Modelling and control of piezoelectric perforated structures}, Proceedings of The Third World Conference On Structural Control. John Wiley, Chichester. F. Casciati : Editor. Vol 3, (2003) 329-336.
\bibitem{Olei} Oleinik O.A., Shamaev G.A., Yosifian G.A. {\it Mathematical problems in elasticity and homogenization}, North  Holland, Amsterdam (1992).
\end{thebibliography}
\end{document}